\documentstyle[11pt,leqno]{article}
\newtheorem{th}{Theorem}[section]

\newtheorem{co}[th]{Corollary}

\title{\Large \bf {Weyl structures with positive Ricci tensor}}
\author{{\sc B. Alexandrov} \hspace{3mm} {\sc S. Ivanov}
\thanks{The authors are supported by Contract MM 809/1998 with the
Ministry of Science and Education of Bulgaria and by Contract 238/1998 with the
University of Sofia "St. Kl. Ohridski".}}
\date{}
\begin{document}
\maketitle
\thispagestyle{empty}
\vspace{2mm}
\begin{center}
Bogdan Alexandrov\\ Sofia University 'St. Kl. Ohridski',\\ Faculty
of Mathematics and Informatics, \\ 5 James Bourchier Blvd, 1164
Sofia, BULGARIA.\\ {\tt alexandrovbt@fmi.uni-sofia.bg} \vspace{4mm}

Stefan Ivanov\\ University of Sofia 'St. Kl. Ohridski',\\ Faculty
of Mathematics and Informatics, \\ 5 James Bourchier Blvd, 1164
Sofia, BULGARIA.\\ {\tt ivanovsp@fmi.uni-sofia.bg}
\end{center}
\vspace{15mm}

\begin{abstract}
We prove the vanishing of the first Betti number on compact
manifolds admitting a Weyl structure whose Ricci tensor satisfies
certain positivity conditions, thus obtaining a Bochner-type
vanishing theorem in Weyl geometry. We also study compact
Hermitian-Weyl manifolds with non-negative symmetric part of the
Ricci tensor of the canonical Weyl connection and show that every
such manifold has first Betti number $b_1 =1$ and Hodge numbers
$h^{p,0} =0$ for  $p>0$, $h^{0,1} =1$, $h^{0,q} =0$ for $q>1$.
\\[15mm]
{\bf Running title:} Weyl structures with positive Ricci tensor
\\[5mm]
\\[5mm]
${\bf MS}$ {\bf classification: } 53C15; 53C55
\end{abstract}
\newpage

\section{Introduction}

A Weyl structure on a conformal manifold $(M,c)$ is a torsion-free
linear connection $\nabla^W$ preserving the conformal structure
$c$. This implies that for any Riemannian metric $g\in c$ there
exists a 1-form $\theta$ such that $\nabla^W g=\theta \otimes g$ .
Conversely, given a 1-form $\theta $ and $g \in c$, there exists a
unique Weyl structure $\nabla ^W$ such that $\nabla^W g=\theta
\otimes g$. An example of a Weyl structure is the Levi-Civita
connection of a metric in $c$. This is the case when $\theta $ is
exact and such Weyl structures are called exact. Similarly,
$\nabla^W$ is called closed if $\theta $ is closed, i.e., if
$\nabla^W$ locally coincides with the Levi-Civita connection of a
metric in $c$. The well known theorem of Gauduchon \cite{G3} says
that if $M$ is compact and at least 3-dimensional, then for every
Weyl structure there exists a unique (up to homothety) metric $g
\in c$, such that the corresponding 1-form $\theta $ is co-closed
with respect to $g$.

Weyl structures arise naturally for example in almost Hermitian
and contact geometry (the connection is determined by the Lee form
and the contact form respectively). More generally, when the
geometry has a preferred 1-form, then the Weyl structure is
determined by the given metric and this 1-form.

Recently there has been considerable interest in Weyl geometry, mainly in
Einstein-Weyl manifolds. A Weyl structure is called Einstein-Weyl if the
symmetric trace-free part of its Ricci tensor vanishes.
The Einstein-Weyl geometry initially
attracted  particular interest in dimension 3 \cite{t3,t4}, but
subsequently
Einstein-Weyl manifolds have been much studied in all dimensions
\cite{t1,t2,PS} (for a nice overview on Einstein-Weyl manifolds see \cite{CP}).

In the present paper we study Weyl structures on compact manifolds, such that
the symmetric part of the Ricci tensor satisfies certain positivity conditions.

In \cite{pisa} Gauduchon showed that a 4-dimensional compact conformal
manifold, which  admits a closed self-dual Weyl structure
with 2-positive normalized Ricci tensor, is conformally equivalent to $S^4$ or
$CP^2$ with their standard conformal structures.
Notice that the 2-positivity
condition of Gauduchon implies the positivity of the Ricci tensor of the Weyl
connection but the converse is not true. In \cite{PS} Pedersen and Swann proved
that on a non-exact compact Einstein-Weyl manifold with vanishing (resp.
non-negative but non-zero) symmetric part of the Ricci tensor the first Betti
number is $b_1 =1$ (resp. $b_1 =0$). Recently the authors showed \cite{AI}
that a compact Hermitian surface with non-negative and positive at some point
symmetric part of the Ricci tensor
of the canonical Weyl connection has $b_1=0$.

In the next theorem we prove that $b_1=0$ for any compact
conformal manifold, which admits a Weyl structure such that the
symmetric part of the Ricci tensor satisfies certain positivity
conditions. Thus we obtain a vanishing theorem of Bochner type in
Weyl geometry.

\begin{th}\label{th1}
Let $(M,c)$ be a compact $n$-dimensional ($n>2$) conformal
manifold, $\nabla ^W$  a non-exact Weyl connection on $(M,c)$
and $g \in c$  the Gauduchon metric of $\nabla ^W$.  Let
$\theta$ be the 1-form determined by $\nabla ^W$ and $g$ and
$\theta^{\#}$  the vector field dual to $\theta$ with respect to
$g$. Suppose the symmetric part
$Ric^W$ of the Ricci tensor of $\nabla^W$ satisfies
\begin{equation}\label{3}
Ric^W (X,X) \ge \frac{(n-2)(n-4)}{8} (|\theta |^2 |X|^2 - \theta (X) ^2)
\end{equation}
for all tangent vectors $X$. Then $b_1  \le 1$. Further,

a) If at some point Inequality (\ref{3})  is strict for $X =
\theta^{\#}$ (i.e. $Ric^W (\theta^{\#} ,\theta^{\#}) >0$), then
$b_1(M) =0$.

b) If $b_1 =1$, then $\theta $ is parallel with respect to the
Levi-Civita connection of $g$ and the universal cover of $(M,g)$
is isometric to ${\bf R} \times N$, where the metric on $N$ has
positive Ricci curvature. In particular, for $n=3$ and $n=4$ the
manifold $N$ is diffeomorphic to $S^{n-1}$.
\end{th}

Note that if $n=4$ Condition (\ref{3}) is reduced to $Ric^W (X,X)
\ge 0$ and if $n=3$ its right-hand side is even non-positive. In
Corollary~\ref{cor3} we also show that any oriented 3-dimensional
manifold, which satisfies the assumptions of Theorem~\ref{th1} and
has $b_1 =1$, is diffeomorphic to $S^1 \times S^2$.

Recall that the canonical Weyl structure on a Hermitian manifold
is determined by the metric and the Lee form. The manifold is
called Hermitian-Weyl if the canonical Weyl structure preserves
the complex structure (see \cite{PPS}). This condition is always
satisfied in the 4-dimensional case, but for higher dimensions it
forces  the manifold to be locally conformally K\"ahler. Every
locally conformally K\"ahler manifold has closed Lee form and a
particular subclass of such manifolds is formed by the generalized
Hopf manifolds, the Hermitian manifolds whose Lee form is parallel
with respect to the Levi-Civita connection. Compact Hermitian
Einstein-Weyl manifolds are classified in \cite{GI} in dimension 4
and the higher dimensional case is studied in detail  in
\cite{PPS}.

We apply Theorem~\ref{th1} to Hermitian Weyl four-manifolds to get

\begin{co}\label{cor1}
Let $(M,J)$ be a compact complex surface. Suppose  there
exists a conformal class of Hermitian metrics such that the
corresponding canonical Weyl structure is non-exact and has
Ricci tensor with everywhere non-negative symmetric part. Then $b_1
\le 1$ and

(i) If $b_1=0$, then (M,J) is a rational surface with  $c_1^2 \ge
0$;

(ii) If $b_1=1$, then $(M,J)$ is a Hopf   surface of class 1.
\end{co}
Note that if the symmetric part of the Ricci tensor of the
canonical Weyl structure is non-negative everywhere and strictly
positive at some point, then the surface is rational with $c_1^2 >
0$ \cite{AI}.

In dimensions greater than four the non-exact Hermitian-Weyl
manifolds with non-negative symmetric part of the Ricci tensor
are generalized Hopf manifolds with respect to the Gauduchon
metric by the results in \cite{V1}. It is well known that the
first Betti number of a generalized Hopf manifold is odd and
therefore these manifolds do not admit any K\"ahler structure. In
general, the first Betti number and the Hodge numbers of a
generalized Hopf manifold are not known  (for surveys on
generalized Hopf manifolds see \cite{V4,DO}). We prove

\begin{th}\label{t01}
Let $(M,g,J)$ be a  2m-dimensional  compact generalized Hopf
manifold and the canonical Hermitian Weyl structure has Ricci
tensor with non-negative symmetric part. Then  $b_1 =1$ and the
Hodge numbers   $h^{p,q}$   satisfy  $$ h^{p,0}=0,
\quad p=1,2,\ldots,m, \qquad h^{0,q}=0, \quad q=2,3,\ldots,m,
\qquad h^{0,1}=1. $$
\end{th}

\section{Ricci-positive Weyl structures in dimension {\boldmath $n$}}

In the following we denote by $<.,.>$ and $|.|$ the pointwise
inner products and norms and by $(.,.)$ and $\| .\|$ -- the global
ones respectively. For the curvature and the Ricci tensor of a
linear connection $\nabla$ we adopt the following convention:
$R(X,Y) = [\nabla _X,\nabla _Y] - \nabla _{[X,Y]},\qquad
\rho(X,Y)= trace \{ Z \longrightarrow R(Z,X)Y \}.$

Let $\nabla^W$ be a Weyl structure on an  $n$-dimensional conformal
manifold $(M,c)$. Let $g \in c$ and $\nabla^W g=\theta \otimes g$.
The connection $\nabla^W$ is given explicitly by
\begin{equation}\label{n1}
\nabla ^W_X Y= \nabla _X Y - \frac{1}{2} \theta  (X)Y
- \frac{1}{2} \theta  (Y)X + \frac{1}{2} g(X,Y)\theta ^{\#},
\end{equation}
where $\nabla $ is the Levi-Civita connection of $g$.
We shall denote by $Ric^W $ the symmetric part of the Ricci tensor
of $\nabla ^W$ and by $k$  the conformal scalar curvature, i.e.
the trace of $Ric ^W$ with respect to $g$. It is easy to obtain
from (\ref{n1}) that the whole Ricci  tensor of $\nabla ^W$ is $Ric^W +
\frac{n}{4} d\theta $,
\begin{equation}\label{1}
Ric^W (X,X) = \rho (X,X) + \frac{n-2}{2} (\nabla _X \theta)(X) -
\frac{n-2}{4} (|\theta |^2 |X|^2 - \theta (X) ^2) -
\frac{1}{2} |X|^2 d^*\theta ,
\end{equation}
\begin{equation}\label{2}
k = s -(n-1)d^* \theta - \frac{(n-1)(n-2)}{4} |\theta |^2,
\end{equation}
where $\rho $ and $s$ are the Ricci tensor and the scalar
curvature of the Levi-Civita connection of $g$ and the  norms and
the codifferential operator are taken with respect to $g$.

\vspace{3mm}
\noindent {\it Proof of Theorem~\ref{th1}:} Let $\widetilde{\nabla }$ be the
connection defined by
\begin{equation}\label{4}
\widetilde{\nabla }_X Y = \nabla _X Y - \frac{n-2}{4} (\theta (X)Y
- g(X,Y)\theta ^{\#}).
\end{equation}
Let $\varphi $ be a 1-form and $\xi$  the dual vector field. It
follows from (\ref{4}) that
\begin{equation}\label{5}
|\widetilde{\nabla } \xi|^2 = |\nabla \xi|^2 -
\frac{n-2}{2}<\nabla _{\theta ^{\#}} \xi,\xi> + \frac{n-2}{2}\theta (\nabla
_{\xi} \xi) + \frac{(n-2)^2}{8} (|\theta |^2 |\xi|^2 - \theta (\xi)^2).
\end{equation}
Since
$$(\nabla _{\theta ^{\#}} \xi,\xi) = \frac{1}{2}(\theta
,d(|\xi|^2)) = \frac{1}{2}\int _M {|\xi|^2 d^* \theta }\, dV =0,$$
by (\ref{1}), (\ref{5}) and the well-known Weitzenb\"ock formula
$$\|d\varphi \|^2 + \|d^* \varphi \|^2 = \|\nabla \xi \|^2 + \int
_M {\rho (\xi,\xi)} \, dV $$ (see e.g. \cite{Bes}, Chapter 1) we
obtain
\begin{equation}\label{6}
\|d\varphi \|^2 + \|d^* \varphi \|^2 = \|\widetilde{\nabla } \xi \|^2 +
\int _M {(Ric^W (\xi,\xi) - \frac{(n-2)(n-4)}{8}(|\theta |^2 |\xi|^2 - \theta
(\xi)^2))} \, dV.
\end{equation}
Now, let $\varphi \not =0$ be a harmonic 1-form. Then (\ref{6})
and (\ref{3}) show that $\widetilde{\nabla }
\xi =0$. Hence, $\nabla \varphi = \frac{n-2}{4} \theta \wedge \varphi$ is
skew-symmetric. But the skew-symmetric part of $\nabla \varphi$ is $d \varphi
=0$. Therefore $\theta = f\varphi $ for some smooth function
$f$ and  $\nabla \varphi = 0$. Thus, the universal cover of
$(M,g)$ is isometric to ${\bf R} \times N$ (see e.g. chapter 4 in
\cite{KN}). We can lift $\theta $, $\varphi $, $\xi$ and $f$ to
${\bf R} \times N$ and we can assume that $|\varphi |=1$, i.e.
$\xi$ is the unit vector field tangent to ${\bf R}$. Since $\theta $ and
$\varphi$ are both co-closed, from $\theta = f\varphi $ we get $df (\xi) =0$.
By (\ref{1}) and (\ref{3}) we have
\begin{equation}\label{9}
\rho (X,X) + \frac{n-2}{2}\varphi (X) df(X) \ge
\frac{(n-2)^2}{8} f^2 (|\varphi |^2 |X|^2 - \varphi (X) ^2).
\end{equation}
Any vector tangent to ${\bf R} \times N$ has the form $X=\lambda
\xi + X^\bot $, where $X^\bot$ is orthogonal to $\xi$, and since
the metric on ${\bf R} \times N$ is the product metric,
$\rho (X,X) = \rho (X^\bot,X^\bot)$. Thus, using (\ref{9}) we
obtain
\begin{equation}\label{10}
\rho (X^\bot ,X^\bot ) + \frac{n-2}{2}\lambda df(X^\bot) \ge
\frac{(n-2)^2}{8} f^2 |X^\bot |^2
\end{equation}
for any $\lambda \in {\bf R}$. Hence,  for any vector $X^\bot $
tangent to $N$ we have $df(X^\bot) =0$. This together with $df(\xi) =0$
means that $f$ is constant and therefore $\theta $ is parallel.
Since the Weyl structure on $M$ is not exact, we have $f \not =0$
(otherwise $\theta =0$). Thus, any harmonic form on $(M,g)$ is a
constant multiple of $\theta $ and so if $b_1 \not =0$, then $b_1
=1$. It follows from (\ref{10}) that the Ricci curvature of $N$ is
positive.
\\
When $n=4$, $N$ is a compact simply connected
3-dimensional manifold admitting a metric of positive Ricci
curvature and by a theorem of Hamilton (cf. \cite{Ham} or
Theorem~5.30 in \cite{Bes}) $N$ is diffeomorphic to $S^3$. When
$n=3$, $N$ is a compact simply connected 2-dimensional manifold
admitting a metric of positive Ricci curvature (i.e. of positive
Gauss curvature) and therefore is diffeomorphic to $S^2$. Thus b)
is proved. \\
If Inequality (\ref{3}) is strict for
$X=\theta^{\#}$ at some point, then it follows from (\ref{6}) that
$\theta$ is not a harmonic form. This means that there are no
harmonic 1-forms, which proves a). \hfill {\bf Q.E.D.}

\begin{co}\label{cor3}
Any oriented 3-dimensional manifold $M$, which satisfies the
assumptions of Theorem~\ref{th1} and has $b_1 \not=0$, is
diffeomorphic  to $S^1 \times S^2$.
\end{co}

\noindent {\it Proof:} By Theorem~\ref{th1} the universal cover of
$M$ is isometric to ${\bf R} \times N$, where $N$ is diffeomorphic
to $S^2$. Any metric on $S^2$ is conformal to the standard one and
therefore orientation-preserving isometries of $N$ are contained
in $SL(2,{\bf C})$. Now, by arguments similar to those in the
proof of Theorem~3.2 in \cite{AGI} we obtain that the fundamental
group of $M$ is isomorphic to ${\bf Z}$ and $M$ is the total space
of a locally trivial fibre bundle over $S^1$ with fibre $S^2$ and
structure group contained in $SL(2,{\bf C})$. Hence, $M$ is
diffeomorphic to $S^1 \times S^2$.   \hfill {\bf Q.E.D.}

\vspace{3mm}
The following result is proved in \cite{PS} for Einstein-Weyl manifolds.
\begin{co}\label{cor2}
If $(M,c)$ is a compact 4-dimensional conformal manifold with a
Weyl structure whose Ricci tensor is non-negative everywhere and
positive at some point symmetric part, then the 1-form $\theta $
given by any choice of a compatible metric must vanish somewhere.
\end{co}

\noindent{\it Proof:} It follows from Theorem~\ref{th1} that $b_1
=0$ and therefore the Euler characteristic of $M$ is positive.
Hence, every 1-form on $M$ vanishes somewhere. \hfill {\bf Q.E.D.}

\section{Ricci-positive Hermitian-Weyl structures}

Let $(M,g,J)$ be a $2m$-dimensional ($m>1$) almost  Hermitian
manifold with almost complex structure $J$ and compatible metric
$g$. Let $\Omega $ be the K\"ahler form of $(M,g,J)$, defined by
$\Omega (X,Y)=g(X,JY)$. Denote by $\theta$ the Lee form of
$(M,g,J)$, $\theta = -\frac{1}{m-1} Jd^* \Omega$ (for a 1-form
$\alpha $, we define $J\alpha = -\alpha\circ J$).

The canonical Weyl structure $\nabla^W$ corresponding to the
almost Hermitian structure $(g,J)$ on $M$ is determined by the
metric $g$ and its Lee form $\theta$. The canonical Weyl structure
does not depend on the choice of the metric in the conformal class
$c$ of $g$, i.e. the  Weyl structure determined by a metric
$\widetilde{g} \in c$ and its Lee form $\widetilde{\theta}$
coincides with $\nabla^W$.

The canonical Weyl structure is called Hermitian-Weyl \cite{PPS}
if it preserves $J$. Since $\nabla^W$ is torsion-free, this shows
that $J$ must be integrable. According to the results in
\cite{Va3}, the canonical Weyl structure on a Hermitian manifold
$(M,g,J)$ is Hermitian-Weyl iff
\begin{equation}\label{200}
d\Omega=\theta\wedge\Omega .
\end{equation}
This condition is always satisfied in the  4-dimensional case. For
higher dimensions it means that the Hermitian manifold $(M,g,J)$
is locally conformally K\"ahler. In particular, $d\theta =0$.
Notice that since the antisymmetric part of the Ricci tensor of
$\nabla^W$ is $\frac{m}{2} d\theta $, on locally conformally
K\"ahler manifolds the Ricci tensor of the canonical Weyl
connection  is symmetric.

\vspace{3mm} \noindent {\it Proof of Corollary~\ref{cor1}.}
Theorem~\ref{th1} implies $b_1 \le 1$. Since the Weyl structure is
not exact, (\ref{2}) shows that the Gauduchon metric has
non-negative but not identically zero scalar curvature. Thus, by
Gauduchon's plurigenera theorem all  plurigenera of $(M,J)$
vanish, cf. Proposition~I.18 in \cite{G4} or \cite{V}. Hence the
Kodaira dimension of $(M,J)$ is $-\infty$ \cite{BPV}.

If $b_1 =0$, then $h^{0,1}=0$ (see e.g. \cite{BPV}). Hence, by the
Castelnuovo criterion (cf. \cite{BPV}) $(M,J)$ must be a rational
surface. The assertion about $c_1^2$ follows by the same arguments
as those in the proof of Theorem~1.3 in \cite{AI}.

If $b_1=1$, then  Theorem~\ref{th1} shows that $(M,J)$ is a
generalized Hopf manifold with respect to the Gauduchon metric.
The recent classification of generalized Hopf surfaces \cite{Bel}
implies that the surface is a Hopf surface of class 1. \hfill {\bf
Q.E.D.}

\vspace{3mm}
For the rest of this section we assume that $(M,g,J)$  is a
Hermitian manifold, such that the canonical Weyl connection
$\nabla^W $ is Hermitian-Weyl.

The Chern connection $\nabla^C$ of $(M,g,J)$ is given by
$$
g(\nabla ^C _X Y,Z) = g(\nabla _X Y,Z) + \frac{1}{2} d\Omega
(JX,Y,Z).$$
This, together with (\ref{n1}) and
(\ref{200}), yields
$$\nabla^C=\nabla^W +\frac{1}{2}\theta\otimes Id
+ \frac{1}{2}J\theta\otimes J,$$
where $Id$ denotes the identity of
$TM$. Consequently,  the curvature tensors $R^C$ and $R^W$ of
$\nabla^C$ and $\nabla^W$ are related by
\begin{equation}\label{tt4}
R^C=R^W + \frac{1}{2} d(J\theta)\otimes J + \frac{1}{2} d\theta \otimes Id .
\end{equation}
Notice that for $m>2$ the last term in (\ref{tt4}) is  absent,
since in this case $d\theta =0$. Let $$R^C(\Omega)
(X,Y)=\frac{1}{2}\sum_{j=1}^{2m} g(R^C (e_i,Je_i)X,Y),$$ where $\{
e_1,\ldots ,e_{2m} \}$ is an orthonormal frame of tangent vectors.
We denote by $k^C$ the symmetric, J-invariant bilinear form
associated to $R^C(\Omega)$ and defined by
$$k^C(X,Y)=R^C(\Omega)(JX,Y).$$

Using  (\ref{tt4}) and the fact that $R^W\circ J=J\circ R^W$, we obtain
\begin{equation}\label{tt5}
k^C=Ric^W+\frac{1}{2}<d(J\theta),\Omega>g.
\end{equation}
This equality is proved in \cite{GI} (Formula (14))  in dimension
4. The proof in higher dimensions is the same, using that $d\theta
=0$.

\vspace{3mm} \noindent {\it Proof of  Theorem~\ref{t01}:} We shall
use Formula (2.8) in \cite{V4}, according to which on every
generalized Hopf manifold
\begin{equation}\label{tt6}
d(J\theta)=|\theta|^2 \Omega + \theta \wedge J\theta.
\end{equation}
Equalities (\ref{tt5}) and (\ref{tt6}) yield
\begin{equation}\label{tt7}
k^C=Ric^W+(m-1)|\theta|^2 g.
\end{equation}
Since $Ric^W \ge 0$ and  $\theta \not =0$,  (\ref{tt7})
implies  $k^C>0$. This allows us to apply the vanishing
theorem for holomorphic $(p,0)$-forms \cite{KW,Ga} and  deduce
that $h^{p,0}=0$, $p=1,2,\ldots,m$.  On the other hand, the Hodge
numbers on a $2m$-dimensional generalized Hopf manifold satisfy
the following relations \cite{Ts}:
\begin{equation}\label{hop}
h^{m,0}=h^{m-1,0}, \quad h^{0,p}=h^{p,0}+h^{p-1,0}, \, p\le
m-1,\quad
%\end{equation}
%\begin{equation}\label{hop1}
2h^{1,0}=b_1 -1, \quad 2h^{0,1}=b_1 +1.
\end{equation}
From $h^{1,0}=0$ and (\ref{hop}) we get  $b_1=1$, $h^{0,1}=1$.
Further, (\ref{hop}) implies  $h^{0,p}=0$ for $p=2,\ldots,m-1$.
Finally, the Kodaira-Serre duality and $h^{m,0}=0$ lead to the
vanishing of $h^{0,m}$, which completes
the proof.
\hfill {\bf
Q.E.D.}

\vspace{3mm} Notice that if $2m=4$ the result of Theorem~\ref{t01}
follows immediately from Corollary~\ref{cor1}.

\vspace{3mm} \noindent {\bf Acknowledgments}  The second author
thanks to The Abdus Salam International Centre for Theoretical
Physics, Trieste, Italy, where the final part of this work was
done. Both authors thank V. Apostolov for his valuable comments
and suggestions.

\end{document}